\documentclass{amsart}
\numberwithin{equation}{section}

\newtheorem{thm}{Theorem}[section]
\newtheorem{lem}[thm]{Lemma}

\begin{document}
\title{Almost-everywhere
convergence of Fourier series for functions in Sobolev spaces}

\author{Ravshan Ashurov}
\address{National University of Uzbekistan named after Mirzo Ulugbek and Institute of Mathematics, Uzbekistan Academy of Science}
\curraddr{Institute of Mathematics, Uzbekistan Academy of Science,
Tashkent, 81 Mirzo Ulugbek str. 100170} \email{ashurovr@gmail.com}

\small

\title[Almost-everywhere convergence]
{Almost-everywhere convergence of Fourier series for functions in
Sobolev spaces}

\begin{abstract}

Let $S_\lambda F(x)$ be the spherical partial sums of the multiple
Fourier series of function $F\in L_2(\mathbb{T}^N)$. We prove
almost-everywhere convergence $S_\lambda F(x)\rightarrow F(x)$ for
functions in Sobolev spaces $H_p^a(\mathbb{T}^N)$ provided $1< p
\leq 2$ and $a> (N-1)(\frac{1}{p}-\frac{1}{2})$. For multiple
Fourier integrals this is well known result of Carbery and Soria
(1988). To prove our result, we first extend the transplantation
technic of Kenig and Tomas (1980) from $L_p$ spaces to $H_p^a$
spaces, then apply it to the Carbery and Soria result.

\end{abstract}+-

\maketitle

\section{Main result.} Let $\{F_n\}$, $n\in \mathbb{Z}^N$, be
the Fourier coefficients of a function $F\in L_2(\mathbb{T}^N)$,
$N\geq2$, i.e.
$$F_n=(2\pi)^{-N}\int\limits_{\mathbb{T}^N}
F(y)e^{-iny}dy,$$ where $\mathbb{T}^N$ is $N$-dimensional torus:
$\mathbb{T}^N = (\pi, \pi]^N$. Consider the spherical partial sums
of the multiple Fourier series:
\begin{equation}\label{SL}
S_\lambda F(x)=\sum\limits_{|n|^2< \lambda}F_n\,e^{inx},
\end{equation}
where $nx=n_1x_1+n_2x_2+...+n_Nx_N$ and
$|n|=\sqrt{n_1^2+n_2^2+...+n_N^2}$.

The aim of this paper is to investigate convergence
almost-everywhere (a.e.) of these partial sums for functions in
Sobolev spaces. To formulate the main result, we need to remind
the definition of Sobolev spaces: the class of functions
$L_p(\mathbb{T}^N)$ which for a given fixed number $a> 0$ make the
norm
\begin{equation}\label{norm}
||F||_{H_p^a(\mathbb{T}^N)}=\big|\big|\sum\limits_{n\in
Z^N}(1+|n|^2)^{\frac{a}{2}}F_n
e^{inx}\big|\big|_{L_p(\mathbb{T}^N)}
\end{equation}
finite is termed the Sobolev class $H_p^a(\mathbb{T}^N)$, $p\geq
1$.

\begin{thm}\label{M} Let $1 < p \leq 2$ and $a>(N-1)(\frac{1}{p}-\frac{1}{2}).$
Then for any function $F\in H^a_p(\mathbb{T}^N)$ the partial sums
$S_\lambda F(x)$ converge to $F(x)$ a.e. in $\mathbb{T}^N$.
Moreover, the maximal operator $S_\star F(x)=
\sup\limits_{\lambda>0}|S_\lambda F(x)|$ has the estimate
\begin{equation}\label{SMax}
||S_\star F||_{L_p(\mathbb{T}^N)}\leq
C_{p,a}||F||_{H^a_p(\mathbb{T}^N)}.
\end{equation}

\end{thm}
This result was previously proved for the multiple Fourier
integrals in \cite{CAR} (see also \cite{AB}). It may be worth
mentioning that in \cite{A} the statement of Theorem \ref{M} was
proved for spectral decompositions of an arbitrary elliptic
pseudo-differential operator, defined in $L_2(\Omega)$,
$\Omega\subseteq \mathbb{R}^N$, of functions in the Sobolev
classes $H_p^a(\Omega)$ with $1 \leq p \leq 2$ and
$a>N(\frac{1}{p}-\frac{1}{2}),$ and moreover, it was shown that
this result is definitive in the class of all elliptic
pseudo-differential operators.

\section{Theorem of C. Kenig and P. Tomas.} An
$L_\infty(\mathbb{R}^N)$ function $M$ is regulated if every point
of $\mathbb{R}^N$ is a Lebesque point of $M$ (see \cite{KP}).
Define for each real number $\lambda>0$ an operator
$E_{M,\lambda}$ on $L_2(\mathbb{R}^N)$ by
$\widehat{E_{M,\lambda}f} (\xi)=M(\xi/\lambda)\hat{f}(\xi)$ and
$S_{M,\lambda}$ on $L_2(\mathbb{T}^N)$ by
$\widehat{S_{M,\lambda}F} (n)=M(n/\lambda)\hat{F}(n)$. The symbols
$E^\star_M f(x)$ and $S^\star_M F(x)$ stand for the corresponding
maximal operators, i.e. $E^\star_M
f(x)=\sup\limits_{\lambda>0}|E_{M,\lambda}f(x)|$ and $S^\star_M
F(x)=\sup\limits_{\lambda>0}|S_{M,\lambda}F(x)|$. A function $M$
is called $p$ - maximal on $\mathbb{R}^N$ (or weak $p$ - maximal
on $\mathbb{R}^N$) if the operator $E^\star_M $ is bounded (or
weakly bounded) on $L_p(\mathbb{R}^N)$; similar for $S^\star_M $
on $L_p(\mathbb{T}^N)$.

\textbf{Theorem} (C. Kenig and P. Tomas \cite{KP}). \emph{Suppose
that $M$ is a regulated $L_\infty (\mathbb{R}^N)$ function and let
$1<p<\infty$. Then $M$ is $p$ - maximal or weak $p$ - maximal on
$\mathbb{R}^N$ if and only if $M$ is $p$ - maximal or weak $p$ -
maximal on $\mathbb{T}^N$}.

\section{Theorem of Carbery and Soria.} Let the symbol
$E_\lambda f(x)$ stands for the spherical partial integrals of
multiple Fourier integrals of a function $f\in L_2(\mathbb{R}^N)$:
\begin{equation}\label{EL} E_\lambda
f(x)=(2\pi)^{-\frac{N}{2}}\int\limits_{|\xi|<\lambda}\hat{f}(\xi)e^{ix\xi}d\xi.
\end{equation}
Here a Fourier transform of function $f$ is defined as
$$
\hat{f}(\xi)=(2\pi)^{-\frac{N}{2}}\int\limits_{\mathbb{R}^N}
f(x)e^{-ix\xi}dx.
$$

Carbery and Soria in their paper \cite{CAR} considered, among the
other problems, almost-everywhere convergence of Fourier integrals
$E_\lambda f(x)$ for functions $f$ in Sobolev spaces
$H_p^a(\mathbb{R}^N)$. Note, the norm in this space has the form
\begin{equation}\label{norm}
||f||_{H_p^a(\mathbb{R}^N)}=\big|\big|\int\limits_{\mathbb{R}^N}(1+|\xi|^2)^{\frac{a}{2}}\hat{f}(\xi)
e^{ix\xi}d\xi\big|\big|_{L_p(\mathbb{R}^N)}.
\end{equation}

\textbf{Theorem} (Carbery and Soria, \cite{CAR}). \emph{Let $N\geq
2$. If $f\in H^a_p(\mathbb{R}^N)$ with}
$$
1<p\leq 2, \quad a>(N-1)(\frac{1}{p}-\frac{1}{2}),
$$
\emph{then $E_\lambda f \rightarrow f$ a.e}.

Our aim is  to obtain this result for Fourier series $S_\lambda
F(x)$ using the idea of the proof of the  theorem of C. Kenig and
P. Tomas. To do this we should  first reduce things to $L_p$
spaces.

Let $\Delta$ be the Laplace operator. Then
$$
E_\lambda
f=(1-\Delta)^{-\frac{a}{2}}E_\lambda(1-\Delta)^{\frac{a}{2}}f=(1-\Delta)^{-\frac{a}{2}}E_\lambda
g, \quad g=(1-\Delta)^{\frac{a}{2}}f\in L_p(\mathbb{R}^N).
$$
Define
$$
E_{(\lambda, a)}g=(1-\Delta)^{-\frac{a}{2}}E_\lambda
g=(2\pi)^{-\frac{N}{2}}\int_{|\xi|^2<\lambda}(1+|\xi|^2)^{-\frac{a}{2}}\hat{g}(\xi)e^{ix\xi}d\xi,
\quad g\in L_p(\mathbb{R}^N).
$$

The operator $S_{(\lambda, a)}G$, $G\in L_p(\mathbb{T}^N)$  on the
torus $\mathbb{T}^N$ is defined similarly:
$$
S_{(\lambda,
a)}G=\sum_{|n|^2<\lambda}(1+|n|^2)^{-\frac{a}{2}}G_ne^{inx}, \quad
G\in L_p(\mathbb{T}^N).
$$

Let $E^\star_{(a)}$ and $S^\star_{(a)}$ be the corresponding
maximal operators.

It is not hard to see, that there is no function $M$ so that we
could write the operators $E_{(\lambda, a)}$ and $S_{(\lambda,
a)}$ in the form of operators $E_{M, \lambda}$ and $S_{M,
\lambda}$. Therefore one can not use here the theorem of C. Kenig
and P. Tomas. Nevertheless, we have the following statement.

\begin{thm}\label{new} Let $1<p<\infty$. If the operator $E^\star_{(a)}$ is bounded (or
weakly bounded) on $L_p(\mathbb{R}^N)$, then the operator
$S^\star_{(a)}$ is also bounded (or weakly bounded) on
$L_p(\mathbb{T}^N)$.
\end{thm}
Note, unlike the theorem of C. Kenig and P. Tomas, this theorem is
not "if and only if" type, but it is just enough for our purpose.

\section{The Lorentz spaces.} To investigate the weak
boundedness it is convenient to introduce the Lorentz space $L_{p,
\infty}(\Omega)$, $\Omega \subseteq \mathbb{R}^N$, $1<p<\infty$,
(or the weak $L_{p}(\Omega)$ space) consisting of all measurable
functions $f$ which make the norm
\begin{equation}\label{Lor}
||f||_{L_{p, \infty}(\Omega)}=\sup\limits_{t>0}
\{t(d_f(t))^{1/p}\}
\end{equation}
finite. In this definition the symbol $d_f(t)$ stands for the
distribution function
$$
d_f(t) = \mu \{x\in \Omega; |f(x)|>t\},
$$
where $\mu (E)$ is the Lebesque measure of the set $E$.

It is not hard to verify, that
$$
t^p\mu \{x\in \Omega; |f(x)|>t\}\leq ||f||^p_{L_{p,
\infty}(\Omega)}\leq ||f||^p_{L_{p}(\Omega)}.
$$

We also need the Lorentz space $L_{p, 1}(\Omega)$ with the norm
$$
||f||_{L_{p, 1}(\Omega)}=p\int\limits_0^\infty
\big(d_f(t)\big)^{1/p}dt.
$$
These both Lorentz spaces are Banach spaces and $(L_{p,
1}(\Omega))^\star=L_{\frac{p}{p-1}, \infty}(\Omega)$ (see
\cite{LG}, p. 52).

Let $1 < p < \infty$ and $q = \frac{p}{p-1}$. Then by virtue of
the Hahn-Banach theorem one has
\begin{equation}\label{I1}
||f||_{L_{p, 1}(\Omega)}=\sup\limits_{||g||_{L_{q,
\infty}(\Omega)\leq 1}} |\int\limits_\Omega fg dx|
\end{equation}
and
\begin{equation}\label{I2}
||f||_{L_{p, \infty}(\Omega)}=\sup\limits_{||g||_{L_{q,
1}(\Omega)\leq 1}} |\int\limits_\Omega fg dx|.
\end{equation}

\section{A linearization.} Define the Banach space $L_p(\Omega,
l^\infty(\mathbb{Z}^+))$ as the collection of all sequences of
$L_p(\Omega)$ functions $\{f_k\}$ such that the norm
$||\sup\limits_k|f_k|||_{L_p(\Omega)}$ is finite. The Banach space
$L_p(\Omega, l^1(\mathbb{Z}^+))$ is defined similarly. Then
$E^\star_{(a)}$ may be viewed as an operator defined on
$L_p(\mathbb{R}^N)$ and taking values in $L_p(\mathbb{R}^N,
l^\infty(\mathbb{Z}^+))$ and $S^\star_{(a)}$ defined on
$L_p(\mathbb{T}^N)$ with values in $L_p(\mathbb{T}^N,
l^\infty(\mathbb{Z}^+))$. Using this reduction and duality, the
following linearization of maximal operators can be proved exactly
by the same way as in \cite{KP}.

\begin{lem}\label{lin} Let $1<p<\infty$ and $q=\frac{p}{p-1}$. The operator $E^\star_{(a)}$ is
bounded in $L_p(\mathbb{R}^N)$ if and only if
\begin{equation}\label{E}
||\sum\limits_k E_{(\lambda_k, a)} f_k||_{L_q(\mathbb{R}^N)}\leq C
||\sum\limits_k|f_k|||_{L_q(\mathbb{R}^N)}, \quad \{f_k\}\in
L_q(\mathbb{R}^N, l^1(\mathbb{Z}^+)),
\end{equation}
uniformly in all sequences of positive reals $\{\lambda_k\}$.

Similarly, the operator $S^\star_{(a)}$ is bounded in
$L_p(\mathbb{T}^N)$ if and only if
\begin{equation}\label{S}
||\sum\limits_k S_{(\lambda_k, a)} F_k||_{L_q(\mathbb{T}^N)}\leq C
||\sum\limits_k|F_k|||_{L_q(\mathbb{T}^N)}, \quad \{F_k\}\in
L_q(\mathbb{T}^N, l^1(\mathbb{Z}^+)),
\end{equation}
uniformly in all sequences of positive reals $\{\lambda_k\}$.
\end{lem}

Similar results hold for weak boundedness of operators
$E^\star_{(a)}$ and $S^\star_{(a)}$. Namely, using the equalities
(\ref{I1}) and (\ref{I2}) one can prove that the pair of
inequalities
\begin{equation}\label{EW}
||E^\star_{(a)}f||_{L_{p, \infty}(\mathbb{R}^N)}\leq
||f||_{L_p(\mathbb{R}^N)}, \quad f\in {L_p(\mathbb{R}^N)},
\end{equation}
\begin{equation}\label{SW}
||S^\star_{(a)}F||_{L_{p, \infty}(\mathbb{T}^N)}\leq
||F||_{L_p(\mathbb{T}^N)}, \quad F\in {L_p(\mathbb{T}^N)},
\end{equation}
is equivalent to the pair of inequalities
\begin{equation}\label{EWN}
||\sum\limits_k E_{(\lambda_k, a)} f_k||_{L_q(\mathbb{R}^N)}\leq C
||\sum\limits_k|f_k|||_{L_{q, 1}(\mathbb{R}^N)}, \quad \{f_k\}\in
L_{q, 1}(\mathbb{R}^N, l^1(\mathbb{Z}^+)),
\end{equation}
\begin{equation}\label{SWN}
||\sum\limits_k S_{(\lambda_k, a)} F_k||_{L_q(\mathbb{T}^N)}\leq C
||\sum\limits_k|F_k|||_{L_{q, 1}(\mathbb{T}^N)}, \quad \{F_k\}\in
L_{q, 1}(\mathbb{T}^N, l^1(\mathbb{Z}^+)).
\end{equation}

\section{Proof of Theorem \ref{new}.} We first prove that
$S^\star_{(a)}$ is bounded in $L_p(\mathbb{T}^N)$ if
$E^\star_{(a)}$ is bounded in $L_p(\mathbb{R}^N)$.

According to the linearization Lemma \ref{lin}, it suffices to
show that from the inequality (\ref{E}) uniformly in all sequences
$\lambda_k$ it follows the inequality (\ref{S}) uniformly in all
sequences $\lambda_k$.

Let $G(x)$ be any continuous periodic function on $\mathbb{T}^N$.
Then
\begin{equation}\label{G}
\lim\limits_{\varepsilon\rightarrow
0}\varepsilon^{\frac{N}{2}}\int\limits_{\mathbb{R}^N}G(x)e^{-\varepsilon\pi
|x|^2} dx= \int\limits_{\mathbb{T}^N} G(x) dx,
\end{equation}
as in Stein and Weiss \cite{ST}, p. 261.

Now suppose that the inequality (\ref{E}) holds true. It suffices
to prove (\ref{S}) for an arbitrary trigonometric polynomials
$F_k(x)=P_k(x)$. Let $Q(x)$ be a trigonometric polynomial on
$\mathbb{T}^N$ and $L_\varepsilon (x)=e^{-\varepsilon\pi |x|^2}$.
The following identity can be proved by virtue of (\ref{G}) (see
Stein and Weiss \cite{ST}, p. 261):
$$
\int\limits_{\mathbb{T}^N} \sum\limits_k S_{(\lambda_k, a)}P_k(x)
Q(x) dx= \lim\limits_{\varepsilon\rightarrow
0}\varepsilon^{\frac{N}{2}}\int\limits_{\mathbb{R}^N}\sum\limits_k
E_{(\lambda_k, a)}(P_k L_{\varepsilon/q}) Q(x)L_{\varepsilon/p}(x)
dx.
$$
In \cite{ST} this identity was proved for the operators $S_{M,
\lambda}$ and $E_{M, \lambda}$. But this proof remains valid if
the multiplier is only bounded and continuous.

Applying the inequality (\ref{E}), one has by virtue of the last
equality
$$
|\int\limits_{\mathbb{T}^N} \sum\limits_k S_{(\lambda_k, a)}P_k(x)
Q(x) dx|\leq C_p \limsup\limits_{\varepsilon\rightarrow
0}\big[\varepsilon^{\frac{N}{2}}||\sum\limits_k
|P_k|L_{\varepsilon/q}||_{L_q(\mathbb{R}^N)}||Q
L_{\varepsilon/p}|||_{L_p(\mathbb{R}^N)}\big]=
$$
$$
=C_p \limsup\limits_{\varepsilon\rightarrow
0}\big[||\varepsilon^{\frac{N}{2}}\sum\limits_k
|P_k|L_{\varepsilon/q}||_{L_q(\mathbb{R}^N)}||\varepsilon^{\frac{N}{2}}Q
L_{\varepsilon/p}|||_{L_p(\mathbb{R}^N)}\big]=
$$
$$
=C ||\sum\limits_k|P_k|||_{L_q(\mathbb{T}^N)}
||Q||_{L_p(\mathbb{T}^N)},
$$
where we used (\ref{G}) in the last equality. Taking the supremum
over all trigonometric polynomials $Q$ with $L_p$ norm 1, we
obtain (\ref{S}), and this completes the proof of the theorem
concerning the boundedness of the operators.

As to the weak boundedness of the operators, the proof is similar:
we again apply the linearization Lemma \ref{lin}, and assuming
(\ref{EWN}), in order to prove (\ref{SWN}) it will suffice to know
that
$$
\lim\limits_{\varepsilon\rightarrow
0}\varepsilon^{\frac{N}{2q}}||\sum\limits_k|P_k|L_{\varepsilon/q}||_{L_{q,
1}(\mathbb{R}^N)} = \sum\limits_k|P_k|||_{L_{q, 1}(\mathbb{T}^N)},
$$
which follows from the Poisson summation formula (see \cite{LG},
p. 225, 231).

\section{Proof of Theorem \ref{M}}. Let $1<p\leq 2$,
$q=\frac{2N}{N-1+2/p}$ and $a=\frac{N-1}{2}.$ Application of the
Marcinkeiwicz interpolation theorem to the estimates, obtained in
\cite{CAR}, paragraphs 3 and 4, gives
$$
||E^\star_{(a)}f||_{L_{q}(\mathbb{R}^N)}\leq C
||f||_{L_p(\mathbb{R}^N)}, \quad f\in {L_p(\mathbb{R}^N)}.
$$
Therefore, $E^\star_{(a)}f$ is finite a.e. on $\mathbb{R}^N$ for
each $f\in {L_p(\mathbb{R}^N)}$, $p>1$. Then by Stein's theorem on
a sequence of translation invariant linear operators: $L_p
\rightarrow L_p$ , $1\leq p \leq 2$ (see, for example, \cite{AAP},
p. 73), one may conclude, that
$$
||E^\star_{(a)}f||_{L_{p, \infty}(\mathbb{R}^N)}\leq C
||f||_{L_p(\mathbb{R}^N)}, \quad f\in {L_p(\mathbb{R}^N)},\quad
p>1,\quad a=\frac{N-1}{2},
$$
By virtue of Theorem \ref{new} we have for the same values of the
parameters $p$ and $a$ an estimate
\begin{equation}\label{S1}
||S^\star_{(a)}F||_{L_{p, \infty}(\mathbb{T}^N)}\leq C
||F||_{L_p(\mathbb{T}^N)}, \quad F\in {L_p(\mathbb{T}^N)}.
\end{equation}

On the other hand, from the estimates, obtained in \cite{CAR},
paragraph 3, it follows
$$
||E^\star_{(a)}f||_{L_{2}(\mathbb{R}^N)}\leq C_a
||f||_{L_2(\mathbb{R}^N)}, \quad f\in {L_2(\mathbb{R}^N)},\quad
a>0,
$$
or by virtue of Theorem \ref{new}
\begin{equation}\label{S2}
||S^\star_{(a)}F||_{L_{2}(\mathbb{T}^N)}\leq C_a
||F||_{L_2(\mathbb{T}^N)}, \quad F\in {L_2(\mathbb{T}^N)},\quad
a>0.
\end{equation}

Now applying first the Marcinkeiwicz interpolation theorem  to the
estimates (\ref{S1}) and (\ref{S2}) (with $a=\frac{N-1}{2}$), we
obtain
$$
||S^\star_{(a)}F||_{L_{p}(\mathbb{T}^N)}\leq C
||F||_{L_p(\mathbb{T}^N)}, \quad F\in {L_p(\mathbb{T}^N)},\quad
p>1,\quad a=\frac{N-1}{2},
$$
then applying to this estimate and (\ref{S2}) (with $a>0$) Stein's
interpolation theorem on an analytic family of linear operators
(see, for example, \cite{AAP}, p. 46) we have
$$
||S^\star_{(a)}F||_{L_{p}(\mathbb{T}^N)}\leq
C_{p,a}||F||_{L_p(\mathbb{R}^N)}, \quad F\in {L_p(\mathbb{T}^N)},
$$
where $1<p\leq 2$ and $a>(N-1)(\frac{1}{p}-\frac{1}{2})$.

Turning back to Sobolev spaces, we rewrite this inequality as
$$
||S_\star G||_{L_p(\mathbb{T}^N)}\leq
C_{p,a}||G||_{H^a_p(\mathbb{T}^N)}, \quad G\in
{H^a_p(\mathbb{T}^N)},
$$
where $p$ and $a$ satisfy the same conditions.

This is the estimate (\ref{SMax}). First part of the statement of
Theorem \ref{M} follows from this estimate by using the standard
technic (see \cite{ST}).

\section{Acknowledgement.} The author conveys thanks to Sh. A.
Alimov for discussions of this result and gratefully acknowledges
S. Umarov (University New Haven, USA) for support and hospitality.

 The author was supported by Foundation for Support of Basic Research of the Republic of Uzbekistan
 (project number is OT-F4-88).

\

\

\bibliographystyle{amsplain}

\end{document}